\documentclass[12pt,reqno,a4wide]{amsart}

\usepackage{tikz} 
\usepackage{calrsfs}
\usepackage[charter]{mathdesign}

\allowdisplaybreaks

\newcommand{\gauss}[2]{\genfrac{[}{]}{0pt}{}{#1}{#2}}

\def\schur{\operatorname{{\sf Schur}}}

\title{Some combinatorial matrices and their LU-decomposition}

\author[H.~Prodinger]{Helmut Prodinger}
\address{H.~Prodinger\\Department of Mathematical Sciences, Mathematics Division\\ Stellenbosch University, Private Bag X1, 7602 Matieland, South Africa}
\email{hproding@sun.ac.za}

\thanks{}

\begin{document}
	\begin{abstract}
		Three combinatorial matrices are considered and their LU-decompositions were found. This is typically done
		by (creative) guessing, and necessary proofs are more or less routine calculations.
	\end{abstract}

	\maketitle

\section{Introduction}

Combinatorial matrices often have beautiful LU-decompositions, which leads also to easy determinant evaluations. It has become a habit of this author to try this decomposition whenever he sees a new such matrix. 

The present paper contains three independent ones collect over the last one or two years.

\section{A matrix from polynomials with bounded roots}

In \cite{KT}  Kirschenhofer and Thuswaldner evaluated the determinant
\begin{equation*}
D_s=\det \Bigl(\frac1{(2l)^2-t^2(2i-1)^2}\Bigr)_{1\le i, l\le s}
\end{equation*}
for $t=1$. Consider the matrix $M$ with entries $1/((2l)^2-t^2(2i-1)^2)$ where $s$ might be a positive integer or infinity. In \cite{KT}, the transposed matrix was considered, but that is immaterial when it comes to the determinant; we will treat the transposed matrix as well, but the results are slightly uglier. 

The aim   is to provide a completely elementary evaluation of this determinant which relies on the LU-decomposition $LU=M$, which is obtained by guessing. The additional parameter $t$ helps with guessing and makes the result even more general. We found these results:
\begin{gather*}
L_{i,j}=\frac{\prod_{k=1}^j\bigl((2j-1)^2t^2-(2k)^2 \bigr)}{\prod_{k=1}^j\bigl((2i-1)^2t^2-(2k)^2 \bigr)}
\frac{(i+j-2)!}{(i-j)!(2j-2)!},\\
U_{j,l}=\frac{t^{2j-2}(-1)^j16^{j-1}(2j-2)!}{\prod_{k=1}^j\bigl((2k-1)^2t^2-(2l)^2 \bigr)
	\prod_{k=1}^{j-1}\bigl((2j-1)^2t^2-(2k)^2 \bigr)}
\frac{(j+l-1)!}{l(l-j)!}.
\end{gather*}

 Note that
\begin{equation*}
\prod_{k=1}^j\bigl((2i-1)^2t^2-(2k)^2 \bigr)=
(-1)^j4^j\frac{\Gamma(j+1-t(i-\tfrac12))}{\Gamma(1-t(i-\tfrac12))}
\frac{\Gamma(j+1+t(i-\tfrac12))}{\Gamma(1+t(i-\tfrac12))}
\end{equation*}
and
\begin{equation*}
\prod_{k=1}^j\bigl((2k-1)^2t^2-(2l)^2 \bigr)=
4^jt^{2j}\frac{\Gamma(j+\tfrac12+\tfrac{l}{t})}{\Gamma(\tfrac12+\tfrac{l}{t})}
\frac{\Gamma(j+\tfrac12-\tfrac{l}{t})}{\Gamma(\tfrac12-\tfrac{l}{t})};
\end{equation*}
using these formul\ae, $L_{i,j}$ resp. $U_{j,l}$ can be written in terms of Gamma functions.

The proof that indeed $\sum_j L_{i,j}U_{j,l}=M_{i,l}$ is within the reach of computer algebra systems (Zeilberger's algorithm). 
An old version of Maple (without extra packages) provides this summation. 

As a bonus, we also state the inverses matrices:
\begin{equation*}
L_{i,j}^{-1}=\frac{\prod\limits_{k=1}^{i-1}((2j-1)^2t^2-(2k)^2)}{\prod\limits_{k=1}^{i-1}((2i-1)^2t^2-(2k)^2)}	
\frac{(-1)^{i+j}(2i-2)!(2j-1)}{(i+j-1)!(i-j)!}
\end{equation*}
and
\begin{equation*}
U_{j,l}^{-1}=\prod_{k=1}^{l-1}\big((2k-1)^2t^2-(2j)^2\big)	
\prod_{k=1}^{l}\big((2l-1)^2t^2-(2k)^2\big)	\frac{(-1)^{j}2j^2}{t^{2l-2}(2l-2)!(j+l)!(l-j)!16^{l-1}};
\end{equation*}
the necessary proofs are again automatic.

Consequently the determinant is
\begin{equation*}
D_s=\prod_{j=1}^s U_{j,j}.
\end{equation*}
For $t=1$, this may be simplified:
\begin{align*}
D_s&=\frac1{s!}\prod_{j=1}^s\frac{(-1)^j16^{j-1}(2j-2)!(2j-1)!}{\prod_{k=1}^j(2k-2j-1)(2k+2j-1)
	\prod_{k=1}^{j-1}(2j-2k-1)(2j+2k-1)}\\
&=\frac1{s!}\prod_{j=1}^s\frac{16^{j-1}(2j-1)!^2}{(4j-1)!!(4j-3)!! }
=\frac{4^s}{s!}\prod_{j=1}^s\frac{32^{j-1}(2j-1)!^4 }{(4j-1)!(4j-2)! }\\
&=\frac{4^s\,16^{s(s-1)}}{s!^2}\bigg/\prod_{j=1}^s\binom{4j}{2j}\binom{4j-2}{2j-1}
=\frac{4^s\,16^{s(s-1)}}{s!^2}\bigg/\prod_{j=1}^{2s}\binom{2j}{j}\\
&=\frac{16^{s(s-1)}}{s!^2}\bigg/\prod_{j=0}^{2s-1}\binom{2j+1}{j};
\end{align*}
the last expression was given in \cite{KT}.
 We used the notation $(2n-1)!!=1\cdot 3\cdot 5\cdots (2n-1)$.
 
 Now we briefly mention the equivalent results for the transposed matrix:
  \begin{gather*}
 	L_{i,j}=\frac{\prod_{k=1}^j\bigl((2k-1)^2t^2-(2j)^2 \bigr)}{\prod_{k=1}^j\bigl((2k-1)^2t^2-(2i)^2 \bigr)}
 	\frac{(i+j-1)!j}{(i-j)!(2j-1)!i},\\
 	U_{j,l}=\frac{t^{2j-2}(-1)^j16^{j-1}(2j-1)!}{\prod_{k=1}^j\bigl((2l-1)^2t^2-(2k)^2 \bigr)
 		\prod_{k=1}^{j-1}\bigl((2k-1)^2t^2-(2j)^2 \bigr)}
 	\frac{(j+l-2)!}{j(l-j)!},
 \end{gather*}

\begin{gather*}
L_{i,j}^{-1}=\frac{\prod_{k=1}^{i-1}((2k-1)^2t^2-(2j)^2)}
{\prod_{k=1}^{i-1}((2k-1)^2t^2-(2i)^2)}
\frac{(-1)^{i+j}(2i)!j^2}{(i-j)!(i+j)!i^2},\\
U_{j,l}^{-1}=\prod_{k=1}^l\big( (2k-1)^2t^2-(2l)^2\big)\prod_{k=1}^{l-1}\big( (2j-1)^2t^2-(2k)^2\big)\\
\times\frac{(2j-1)!l!(-1)^j}{t^{2l-2}16^{l-1}(2l-1)!(l+j-1)!(l-j)!(l-1)!}.
 \end{gather*}

\section{Lehmer's tridiagonal matrix}

Ekhad and Zeilberger \cite{EZ} have unearthed Lehmer's \cite{Lehmer} tridiagonal $n\times n$ matrix $M=M(n)$ with entries
\begin{equation*}
	M_{i,j}=\begin{cases}
		1 & \text{if } i=j,\\
		z^{1/2}q^{(i-1)/2}& \text{if } i=j-1,\\
		z^{1/2}q^{(i-2)/2}& \text{if } i=j+1,\\
		0&\text{otherwise.}
	\end{cases}
\end{equation*}

Note the similarity to Schur's determinant
\begin{equation*}   
	\schur(x):=\left|
	\begin{matrix}   
		1&xq^{1+m}&&&&\dots\\
		-1&1&xq^{2+m}&&&\dots\\
		&-1&1&xq^{3+m}&&\dots\\
		&&-1&1&xq^{4+m}&\dots\\
		&&&\ddots&\ddots&\ddots
	\end{matrix}
	\right|
\end{equation*}
that was used to great success in \cite{IPS}. This success was based on the two recursions
\begin{equation*}   
\schur(x)=\schur(xq)+xq^{1+m}\schur(xq^2)
\end{equation*}
and, with
\begin{equation*}   
\schur(x)=\sum_{n\ge0}a_nx^n, 
\end{equation*}
by
\begin{equation*}   
a_n=q^na_n+q^{1+m}q^{2n-2}a_{n-1},
\end{equation*}
 leading to
\begin{equation*}   
a_n=\frac{q^{n^2+mn}}{(1-q)(1-q^2)\dots(1-q^n)}.
\end{equation*} 
Schur's (and thus Lehmer's) determinant plays an instrumental part in proving the celebrated Rogers-Ramanujan identities and generalizations.

Lehmer \cite{Lehmer} has computed the limit for $n\to \infty$ of the determinant of the matrix $M(n)$. Ekhad and Zeilberger \cite{EZ} have generalized this result by computing the determinant of the finite matrix $M(n)$. Furthermore, a lively account of how modern computer algebra leads to a solution   was given. Most prominently, the celebrated $q$-Zeilberger algorithm \cite{PWZ} and creative guessing were used.

In this section, the determinant in question is obtained by computing the LU-decomposition $LU=M$. This is done with a computer, and the exact form of $L$ and $U$ is obtained by guessing. A proof that this is indeed the LU-decomposition is then a routine calculation. From it, the determinant in question is computed by multiplying the diagonal  elements of the matrix $U$. By telescoping, the final result is then quite attractive, as already stated and proved by Ekhad and Zeilberger \cite{EZ}.

We use standard notation~\cite{Andrews76}: $(x;q)_n=(1-x)(1-xq)\dots(1-xq^{n-1})$, and the Gaussian $q$-binomial coefficients
${\gauss nk}=\frac{(q;q)_{n}}{(q;q)_{k}(q;q)_{n-k}}$.

\subsection{The LU-decomposition of $M$}

Let
\begin{equation*}
	\lambda(j):=\sum_{0\le k\le j/2}\gauss{j-k}{k}(-1)^kq^{k(k-1)}z^k.
\end{equation*}
It follows from the basic recursion of the Gaussian $q$-binomial coefficients \cite{Andrews76} that
\begin{equation}\label{recu}
\lambda(j)=\lambda(j-1)-zq^{j-2}\lambda(j-2).
\end{equation}

Then we have 
\begin{equation*}
	U_{j,j}=\frac{\lambda(j)}{\lambda(j-1)},\quad\text{\quad}
	U_{j,j+1}=z^{1/2}q^{(j-1)/2},
\end{equation*}
and all other entries in the $U$-matrix are zero. Further,
\begin{equation*}
	L_{j,j}=1,\quad\text{\quad}
	L_{j+1,j}=z^{1/2}q^{(j-1)/2}\frac{\lambda(j-1)}{\lambda(j)},
\end{equation*}
and all other entries in the $L$-matrix are zero.

The typical element of the product $(LU)_{i,j}$, that is
\begin{equation*}
	\sum_{1\le k\le n}L_{i,k}U_{k,j}
\end{equation*}
is almost always zero; the exceptions are as follows: If $i=j$,  then we get
\begin{equation*}
	L_{j,j}U_{j,j}+L_{j,j-1}U_{j-1,j}=\frac{\lambda(j)+zq^{j-2}\lambda(j-2)}{\lambda(j-1)}=1,
\end{equation*}
because of the above recursion (\ref{recu}). If $i=j-1$,  then we get
\begin{equation*}
	L_{j-1,j-1}U_{j-1,j}+L_{j-1,j-2}U_{j-2,j}=z^{1/2}q^{(j-2)/2},
\end{equation*}
and if $i=j+1$,  then we get
\begin{equation*}
	L_{j+1,j+1}U_{j+1,j}+L_{j+1,j}U_{j,j}=z^{1/2}q^{(j-1)/2}\frac{\lambda(j-1)}{\lambda(j)}\frac{\lambda(j)}{\lambda(j-1)}
	=z^{1/2}q^{(j-1)/2}.
\end{equation*}
This proves that indeed $LU=M$. Therefore  for the determinant of the Lehmer matrix $M$
we obtain the expression
\begin{equation*}
	\prod_{j=1}^n\frac{\lambda(j)}{\lambda(j-1)}=\frac{\lambda(n)}{\lambda(0)}=
	\sum_{0\le k\le n/2}\gauss{n-k}{k}(-1)^kq^{k(k-1)}z^k.
\end{equation*}
Taking the limit $n\to\infty$,  leads
to the old result by Lehmer for the determinant of the infinite matrix:
\begin{equation*}
	\lim_{n\to\infty}\det(M(n))=\sum_{k\ge0}\frac{(-1)^kq^{k(k-1)}z^k}{(q;q)_k}.
\end{equation*}

\textbf{Remarks.}

1. For $q=1$, Lehmer's determinant plays a role when enumerating   lattice paths (Dyck paths) of bounded height, or 
planar trees of bounded height, see \cite{deBrKnRi72, Knuth-Selected, HHPW}.

2. Recursions as in (\ref{recu}) have been studied in \cite{Andrews-fibo, Cigler03, PaPr03} and are linked to so-called
Schur polynomials \cite{Schur17}.
 
 \section{Matrices for Fibonacci polynomials}
 
 Cigler~\cite{Cigler} introduced several matrices that have Fibonacci polynomials as determinants; we will only treat two of them
 as showcases. 
 
 The Fibonacci polynomials are
 \begin{equation*}
F_n(x)=\sum_h\binom{n-h}{h}x^{2n-k};
 \end{equation*}
 our answers will come out in terms of the related polynomials
 \begin{equation*}
f_n=\sum_h \binom{n+h}{2h}X^{h}
 \end{equation*}
 where we write $X=x^2$ for simplicity. It is easy to check that
 \begin{equation*}
 	f_n=(X+2)f_{n-1}-f_{n-2},
 \end{equation*}
 for instance by comparing coefficients.
 
 The first matrix is
 \begin{equation*}
M=\Bigg(\binom{i-1}{j}X+\binom{i+1}{j+1}\Bigg)_{0\le i,j<n}
 \end{equation*}
 and we will determine its LU-decomposition $M=LU$.
 
 We obtained
 \begin{equation*}
L_{i,j}=\frac{\binom{i+1}{j+1}\sum_h\binom{j+h}{2h}X^h+\binom{i}{j}\sum_h\binom{j+h}{2h-1}X^h}{\sum_h\binom{j+1+h}{2h}X^h}
 \end{equation*}
 and
 \begin{gather*}
U_{j,j}=\frac{\sum_h\binom{j+1+h}{2h}X^h}{\sum_h\binom{j+h}{2h}X^h},\\
 	U_{j,l}=(-1)^{j+l}\frac{\sum_h\binom{j+h}{2h-1}X^h}{\sum_h\binom{j+h}{2h}X^h}  	,\quad j<l.
 \end{gather*}

For a proof, we do this computation:
\begin{align*}
\sum_{j}L_{i,j}U_{j,l}&=
\binom{i}{l+1}+\binom il \frac{f_{l+1}}{f_l}
+\sum_{0\le j<l}(-1)^{j+l}\binom{i}{j+1}\\&
+\sum_{0\le j<l}\binom ij \frac{f_{j+1}}{f_j}
-\sum_{0\le j<l}\binom{i}{j+1}\frac{f_j}{f_{j+1}}
-\sum_{0\le j<l}\binom ij\\
&=\binom{i}{l+1}+\binom il \frac{f_{l+1}}{f_l}+\binom il \frac{f_{l-1}}{f_l}-\binom il
+X\sum_{0\le j<l}(-1)^{j+l}\binom ij\\
&=\binom{i}{l+1}+(X+2)\binom il-\binom il-X\binom {i-1}{l-1}\\
&=\binom{i+1}{l+1}+X\binom {i-1}l.
\end{align*}

The determinant is then $U_{0,0}U_{1,1}\dots U_{n-1,n-1}$, and by telescoping
\begin{equation*}
\sum_h\binom{n+h}{2h}X^h=\sum_h\binom{2n-h}{h}x^{2n-2h}=F_{2n}(x).
\end{equation*}
 
For completeness, we also factor the transposed matrix as $LU=M^t$:

\begin{gather*}
L_{i,j}=(-1)^{i+j}\frac{\sum_h\binom{j+h}{2h-1}X^h}{\sum_h\binom{j+1+h}{2h-1}X^h},\quad\text{for } j<i,\\
L_{i,i}=1,
\end{gather*}
and
\begin{equation*}
U_{j,l}=\frac{\binom lj \sum_h \binom{j+h}{2h-1}X^h
	+\binom {l+1}{j+1} \sum_h \binom{j+h}{2h}X^h
	}{\sum_h\binom{j+h}{2h}X^h}.
\end{equation*}

 Now we move to the second matrix:
  \begin{equation*}
 	M=\Bigg(\binom{i}{j}X+\binom{i+2}{j+1}\Bigg)_{0\le i,j<n}.
 \end{equation*}
 
 We find
 \begin{equation*}
L_{i,j}=\frac{\binom{i+1}{j+1}\sum_h\binom{j+1+h}{2h+1}X^h
	+\binom{i}{j}\sum_h\binom{j+1+h}{2h}X^h	}
{\sum_h \binom{j+2+h}{2h+1}X^h}
 \end{equation*}
 and
 \begin{gather*}
U_{j,j}=\frac{\sum_{h}\binom{j+2+h}{2h+1}X^h}{\sum_{h}\binom{j+1+h}{2h+1}X^h},\\
U_{j,j+1}=1, \qquad U_{j,l}=0\quad  \text{for } l\ge j+2.
 	\end{gather*}
 
 For a proof, we compute
 \begin{align*}
\sum_jL_{i,j}U_{j,l}&=\frac{\binom{i+1}{l+1}\sum_h\binom{l+1+h}{2h+1}X^h
	+\binom{i}{l}\sum_h\binom{l+1+h}{2h}X^h	}
{\sum_h \binom{l+1+h}{2h+1}X^h}\\
&+\frac{\binom{i+1}{l}\sum_h\binom{l+h}{2h+1}X^h
	+\binom{i}{l-1}\sum_h\binom{l+h}{2h}X^h	}
{\sum_h \binom{l+1+h}{2h+1}X^h}
 \end{align*}
 and
 \begin{align*}
 	\sum_h \binom{l+1+h}{2h+1}X^h&\sum_jL_{i,j}U_{j,l}=\binom{i+2}{l+1}\sum_h\binom{l+1+h}{2h+1}X^h-\binom{i+1}{l}\sum_h\binom{l+1+h}{2h+1}X^h\\*
 	& 		+\binom{i}{l}\sum_h\binom{l+1+h}{2h}X^h	+\binom{i+1}{l}\sum_h\binom{l+h}{2h+1}X^h\\
 	& 		+\binom{i+1}{l}\sum_h\binom{l+h}{2h}X^h	 		-\binom{i}{l}\sum_h\binom{l+h}{2h}X^h	\\
&=\binom{i+2}{l+1}\sum_h\binom{l+1+h}{2h+1}X^h 		+\binom{i}{l}\sum_h\binom{l+h}{2h-1}X^h	\\*
&=\binom{i+2}{l+1}\sum_h\binom{l+1+h}{2h+1}X^h 		+\binom{i}{l}X\sum_h\binom{l+1+h}{2h+1}X^h	\\*
	 \end{align*}
	  and therefore
 \begin{align*}
 	\sum_jL_{i,j}U_{j,l}
 	&=\binom{i+2}{l+1}		+\binom{i}{l}X,
 	\\ \end{align*}
 as required. The determinant is then
 \begin{equation*}
\sum_{h}\binom{n+1+h}{2h+1}X^h=
\sum_{h}\binom{n+1+h}{n-h}X^h=
\sum_{j}\binom{2n+1+h}{j}x^{2n-2j}=x^{-1}F_{2n+1}(x^2).
 \end{equation*}
 
 For the transposed matrix $LU=M^t$, we find
 \begin{gather*}
 L_{i,i-1}=\frac{\sum_h\binom{i+h}{2h+1}X^h}{\sum_h\binom{i+1+h}{2h+1}X^h},\\
 L_{i,i}=1,\qquad  L_{i,j}=0\quad\text{for } j<i-1,
 \end{gather*}
 and
 \begin{equation*}
U_{j,l}=\frac{\binom{l+1}{j+1}\sum_h\binom{j+1+h}{2h+1}X^h+\binom{l}{j}\sum_h\binom{j+1+h}{2h}X^h}{\sum_h\binom{j+1+h}{2h+1}X^h}.
 \end{equation*}

\end{document}